\documentclass[11pt]{article}
\usepackage{a4wide,amsmath,amssymb,array,xspace,enumerate,bm,graphicx,xspace}
\begin{document}
\newcommand{\AW}{\textbf{AW}\xspace}
 \newcommand{\tf}[2]{\tfrac{#1}{#2}}
\newcommand{\ST}{\texttt{ST}\xspace}
\newcommand{\tr}{\mathop{\text{trace}}}
\newcommand{\minimize}{\mathop{\text{minimize}}}
\newcommand{\maximize}{\mathop{\text{maximize}}}
\newtheorem{theorem}{Theorem}
\newtheorem{lemma}{Lemma}
\newtheorem{conjecture}{Conjecture}
\newtheorem{corollary}{Corollary}
\newcounter{remark}
\newcommand{\Real}{\mathop{\text{Re}}}
\newcommand{\remark}{\stepcounter{remark}\paragraph{Remark \theremark.}}
\newcommand{\bb}{2.4449779}
\newcommand{\bP}{\bar P}
\newcommand{\bV}{\bar V}
\newcommand{\Imag}{\mathop{\text{Im}}}
\newcommand{\diag}{\mathop{\text{diag}}}
\stepcounter{footnote}
\author{Richard Weber\footnote{Statistical Laboratory, Centre for
    Mathematical Sciences, Wilberforce Road, Cambridge CB2 0WB, rrw1
    at cam.ac.uk. This paper is a slightly updated version of a 
    paper that has been on
    the author's website since 23 November 2006.}}
\title{\LARGE \bfseries The optimal strategy for symmetric rendezvous
  search on $K_3$}
\date{}
\newpage
\maketitle

\begin{abstract}
In the symmetric rendezvous search game played on $K_n$ (the
completely connected graph on $n$ vertices) two players are initially
placed at two distinct vertices (called locations). The game is played
in discrete steps and at each step each player can either stay where
he is or move to a different location. The players share no common
labelling of the locations. They wish to minimize the expected number
of steps until they first meet. Rendezvous search games of this type
were first proposed by Steve Alpern in 1976. They are simple to
describe, and have received considerable attention in the popular
press as they model problems that are familiar in real life. They are
notoriously difficult to analyse.
Our solution of the symmetric rendezvous game on $K_3$ makes this the
first interesting game of its type to be solved. It establishes the 20
year old conjecture that the Anderson-Weber strategy is optimal.
\medskip


\end{abstract}

\section{Symmetric rendezvous search on $K_3$}
In the symmetric rendezvous search game played on $K_n$ (the
completely connected graph on $n$ vertices) two players are initially
placed at two distinct vertices (called locations). The game is played
in discrete steps, and at each step each player can either stay where
he is or move to another location. The players wish to meet as quickly
as possible.  They are required to use an identical strategy, which
must involve some random moves or else the players will never
meet. They have no common labelling of the locations, so a given
player must choose the probabilities with which he moves to each of
the locations at step $k$ as only a function of where he has been at
previous steps.

Let $T$, $w$ and $w_k$ denote respectively the number of the step on
which the players meet, the minimum achievable value of $ET$, and the
minimum achievable value of $E[\min\{T,k+1\}]=\sum_{i=0}^kP(T> i)$. We
call $w$ the `rendezvous value' of the game. A long-standing conjecture
of Anderson and Weber (1990) is that for symmetric rendezvous search
on $K_3$ the rendezvous value is $w=\tf 52$. This rendezvous value is
achieved by a type of strategy which is now commonly known as the
Anderson--Weber strategy (\AW). For rendezvous search on $K_n$ the \AW
strategy specifies that in blocks of $n-1$ consecutive steps the
players should randomize between staying at their initial location and
touring the other $n-1$ locations in random order.  On $K_3$ this
means that in each successive block of two steps, each player should,
independently of the other, either stay at his initial location or
tour the other two locations in random order, doing these with
respective probabilities $\tf 13$ and $\tf 23$. The rendezvous value
with this strategy is $ET=\tf 52$.

Rendezvous search problems have a long history. One finds such a
problem in the `Quo Vadis' problem of Mosteller (1965) and recently as
`Aisle Miles' (O'Hare, 2006). The first formal presentation of our
problem is due to Alpern (1976), who states it as his `Telephone
Problem'.  ``Imagine that in each of two rooms, there are $n$
telephones randomly strewn about. They are connected in a pairwise
fashion by $n$ wires. At discrete times $t=0,1,\ldots\,$, players in
each room pick up a phone and say `hello'. They wish to minimize the
time $t$ when they first pick up paired phones and can communicate.''
The \AW strategy was conjectured to be optimal for $K_3$ by Anderson
and Weber (1990), who proved its optimality for $K_2$.  Subsequently,
there have been proofs that \AW is optimal for $K_3$ within restricted
classes of Markovian strategies, such as those that must repeat in
successive blocks of $k$ steps, where $k$ is small. See Alpern and
Pikounis (2000) (for optimality of \AW amongst $2$-Markovian
strategies for rendezvous on $K_3$), and Fan (2009) (for optimality of
\AW amongst $4$-Markovian strategies for rendezvous on $K_3$, and
amongst $3$-Markovian strategies for rendezvous on $K_4$).

The rest of the paper concerns symmetric rendezvous search on
$K_3$. In Section 2 we prove that \AW is optimal. The symmetric
rendezvous game on $K_3$ becomes the first interesting game of its
type to be fully solved.  In Section 3 we discuss the thinking that
led to discovery of this proof. Section 4 discusses some
generalizations and intriguing open problems.

\section{Optimality of the Anderson--Weber strategy}

Recall that $T$ denotes the step on which the players meet.
Let us begin by commenting that the \AW strategy does not minimize
$P(T>i)$ for all $i=1,2,\ldots\,$. In particular, the \AW strategy
produces $P(T> 4)=1/9$. However, one can find a strategy such that
$P(T> 4)=1/10$.  This is somewhat of a surprise and shows that
$ET=\sum_{i=0}^\infty P(T> i)$ cannot be minimized by minimizing each
term of the sum simultaneously.

With Junjie Fan, we have gained greater computational experience of
the problem and have been motivated to make the conjecture that \AW
achieves $w_k$ for all $k$, i.e., minimizes the truncated sum
$\sum_{i=0}^k P(T> i)$.  In the following section we prove this is
true.  Theorem~1 states that $\{w_k\}_{k=0}^\infty=\{1,\tf
53,2,\tf {20}{9},\tf 73,\tf {65}{27},\ldots\}$ with $w_k\rightarrow
\tf 52$.

\begin{theorem} The Anderson--Weber strategy is optimal for the
symmetric rendezvous search game on $K_3$, minimizing
$E[\min\{T,k+1\}]$ to $w_k$ for all $k=1,2,\ldots\,$, where
\begin{equation}
w_k =\left\{\begin{array}{ll}
\tf 52 -\tf 52 3^{-\tf{k+1}{2}}\,, &\quad \text{when $k$ is odd,}\\[5pt]
\tf 52 -\tf 32 3^{-\tf k2}\,, &\quad \text{when $k$ is even.}
	   \end{array}
\right.\label{wkcalc}
\end{equation}
Consequently, the minimal achievable value of $ET$ is $w=\tf 52$.
\end{theorem}

\paragraph{Proof}
Throughout most of the following a subscript $k$ on a vector means
that its length is $3^k$. A subscript $k$ on a matrix means that it is
$3^k\times 3^k$.  Let
\[
B_k=B_1\otimes B_{k-1}\,,\quad\text{where }B_1=\begin{pmatrix}
1 & 1 & 0\\
0 & 1 & 1\\
1 & 0 & 1
\end{pmatrix}\,.
\]
Here `$\otimes$' denotes the Kronecker product.  We label the rows and
columns of $B_1$ as $0,1,2$. Suppose the locations are arranged on a
circle and that players have a common notion of
clockwise.\footnote{Readers familiar with the problem will be aware
  that it might make a difference whether or not the players are
  equipped with a common notion of clockwise. We assume that they
  are. However, since we show that the \AW strategy cannot be bettered
  and this strategy makes no use of the clockwise information, \AW is
  also optimal if the players do not have a common notion of
  clockwise.}  Suppose player II is initially placed one position
clockwise of player I. Then $B_1(i,j)$ is an indicator for the event
that they do not meet when at the first step player I moves $i$
positions clockwise from his initial location, and player II moves $j$
positions clockwise from his initial location. $B^\top$ contains the
indicators for the same event, but when player II starts two positions
clockwise of player I.  Let $1_k$ denote the length $3^k$ vector of
1s.  Since the starting position of player II is randomly chosen, the
problem of minimizing the probability of not having met after the
first step is that of minimizing
\[
p^\top\bigl(\tf 12(B_1+B_1^\top)\bigr)p\,,
\]
over $p\in\Delta_1$, where
$\Delta_k=\{p\,:\,p\in\mathbb{R}^{3k},\,p\geq 0 \text{ and } 1_k^\top
p=1\}$.  Similarly, the $9$ rows and $9$ columns of $B_2$ can be
labelled as $0,\ldots,8$ (base 10), and also
$00,01,02,10,11,12,20,21,22$ (base 3). The base 3 labelling is
helpful, for we may understand $B_2(i_1i_2,j_1j_2)$ as an indicator
for the event that the players do not meet when at his first and
second steps player I moves to locations that are respectively $i_1$
and $i_2$ positions clockwise from his initial position, and player II
moves to locations that are respectively $j_1$ and $j_2$ positions
clockwise from his initial position. The problem of minimizing the
probability that they have not met after $k$ steps is that of
minimizing
\[
p^\top\bigl(\tf 12(B_k+B_k^\top)\bigr)p\,.
\]

In this manner we can also formulate the problem of minimizing
$E[\min\{T,k+1\}]$. Let $J_k$ be the $3^k\times 3^k$ matrix that is
all 1s and let
\begin{align} 
M_1&=J_1+B_1\nonumber\\ M_k&=J_k+B_1\otimes M_{k-1}\nonumber\\ &=
J_k+B_1\otimes J_{k-1}+\cdots+B_{k-1}\otimes J_1+B_k\,.
\end{align}
Then
\[
w_k=\min_{p\in\Delta_k}\Bigl\{p^\top M_kp\Bigr\}
=\min_{p\in\Delta_k}\Bigl\{ \tf 12 p^\top(M_k+M_k^\top)p\Bigr\}\,.
\]

It is a difficult problem to find the minimizing $p$, because $\tf
12(M_k+M_k^\top)$ is not positive semidefinite once $k\geq 2$. The
quadratic form $p^\top \bigl(\tf 12(M_k+M_k^\top)\bigr)p$ has many
local minima that are not global minimums. For example, the strategy
which randomizes equally over the 3 locations at each step, taking
$p^\top=1_3^\top/3^k=(1,1\ldots,1)/3^k$, is a local minimum of this quadratic form.

Consider, for example, $k=2$. To show that $w_2=2$ we must minimize
$p^\top(M_2+M_2^\top)p$. However, the eigenvalues of $\tf
12(M_2+M_2^\top)$ are $\{19,\tf{5}{2},\tf{5}{2},1,1, 1,1,-\tf{1}{2} ,-
\tf{1}{2}\}$, so this matrix is not positive semidefinite. In general,
the minimization over $x$ of a quadratic form such as $x^\top Ax$ is
$\cal{NP}$-hard if $A$ is not positive semidefinite. An alternative
approach might be to try to show that $\tf 12(M_2+M_2^\top)-2J_2$ is a
copositive matrix. For general $k$ , we would wish to show that
$x^\top\bigl(\tf 12(M_k+M_k^\top)-w_kJ_k\bigr)x\geq 0$ for all $x\geq
0$, where $\{w_k\}_{k=1}^\infty=\{\tf 53,2,\tf {20}{9},\tf 73,\tf
{65}{27},\ldots\}$ are the values obtained by the Anderson--Weber
strategy.  However, to check copositivity numerically is also
$\cal{NP}$-hard.

The key idea in this proof is to exhibit a matrix $H_k$ such that
$M_k\geq H_k\geq 0$, where $H_k$ is positive semidefinite (denoted
$H_k\succeq 0$) and $p^\top H_kp$ is minimized over $p\in\Delta_k$ to
$w_k$.  Since $p$ is nonnegative we must have $p^\top M_kp\geq p^\top
H_kp\geq w_k$ for all $p$. For example, we may take \small
\[
\text{\normalsize $M_2=$}\begin{pmatrix}
3 & 3 & 2 & 3 & 3 & 2 & 1 & 
   1 & 1 \\[2pt] 2 & 3 & 3 & 2 & 3 & 
   3 & 1 & 1 & 1 \\[2pt] 3 & 2 & 3 & 
   3 & 2 & 3 & 1 & 1 & 1 \\[2pt] 1 & 
   1 & 1 & 3 & 3 & 2 & 3 & 3 & 
   2 \\[2pt] 1 & 1 & 1 & 2 & 3 & 3 & 
   2 & 3 & 3 \\[2pt] 1 & 1 & 1 & 3 & 
   2 & 3 & 3 & 2 & 3 \\[2pt] 3 & 3 & 
   2 & 1 & 1 & 1 & 3 & 3 & 2 \\[2pt] 
   2 & 3 & 3 & 1 & 1 & 1 & 2 & 
   3 & 3 \\[2pt] 3 & 2 & 3 & 1 & 1 & 
   1 & 3 & 2 & 3
    \end{pmatrix}
\text{\normalsize $\geq H_2=$}
\begin{pmatrix}
3 & 3 & 2 & 3 & 3 & 2 & 1 & 
   1 & 0 \\[2pt] 2 & 3 & 3 & 2 & 3 & 
   3 & 0 & 1 & 1 \\[2pt] 3 & 2 & 3 & 
   3 & 2 & 3 & 1 & 0 & 1 \\[2pt] 1 & 
   1 & 0 & 3 & 3 & 2 & 3 & 3 & 
   2 \\[2pt] 0 & 1 & 1 & 2 & 3 & 3 & 
   2 & 3 & 3 \\[2pt] 1 & 0 & 1 & 3 & 
   2 & 3 & 3 & 2 & 3 \\[2pt] 3 & 3 & 
   2 & 1 & 1 & 0 & 3 & 3 & 2 \\[2pt] 
   2 & 3 & 3 & 0 & 1 & 1 & 2 & 
   3 & 3 \\[2pt] 3 & 2 & 3 & 1 & 0 & 
   1 & 3 & 2 & 3
\end{pmatrix}.
\]\normalsize
where $\tf 12(H_2+H_2^\top)$ is positive semidefinite, with
eigenvalues $\{18, 3, 3, \tf 32, \tf 32, 0, 0, 0, 0\}$. The minimum
value of $p^\top H_2p$ is 2.
\smallskip

We restrict our search for $H_k$ to matrices of a special form.  For
$i=0,\ldots,3^k-1$ we write $i_{\text{base}\,3}=i_1\cdots\,i_k$
(always keeping $k$ digits, including leading 0s when $i\leq
3^{k-1}-1$); so $i_1,\ldots,i_k\in\{0,1,2\}$. We define
\[
P_i=P_{i_1\cdots\,i_k}=P_1^{i_1}\otimes\cdots\otimes P_1^{i_k}\,,
\]
where
\[
P_1=\begin{pmatrix}0&1&0\\0&0&1\\1&0&0\end{pmatrix}\,.
\]
Note that the subscript is now used for something other than the size
of the matrix. It will always be easy for the reader to know the $k$
for which $P_i$ is $3^k\times 3^k$ by context.  Observe that
$M_k=\sum_i m_k(i)P_i$, where $m_k$ is the first row of $M_k$. This
motivates a search for an appropriate $H_k$ amongst those of the form
\[
H_k = \sum_{i=0}^{3^k-1} x_k(i) P_i\,.
\]
The condition $M_k\geq H_k$ is equivalent to $m_k\geq x_k$. In the
example above, $H_2=\sum_i x_2(i)P_i$, where
$x_2=(3,3,2,3,3,2,1,1,0)$, the first row of $H_2$.

Let us observe that the matrices $P_0,\ldots,P_{3^k-1}$ commute with
one another and so have a common set of eigenvectors.  Also, $P_i^\top
=P_{i'}$, where $i'_{\text{base}\,3}=i'_1\cdots\, i'_k$ is obtained from
$i_{\text{base}\,3}=i_1\cdots\, i_k$ by letting $i'_j$ be $0,2,1$ as
$i_j$ is $0,1,2$, respectively. 

Let the columns of the matrices $U_k$ and $W_k$ contain the common
eigenvectors of the $\tf 12(P_i+P_i^\top)$ and also of $\tf
12(M_k+M_k^\top)$. The columns of $W_k$ are eigenvectors with
eigenvalues of $0$. We shall now argue that the condition $\tf
12(H_k+H_k^\top)\succeq 0$ is equivalent to $U_kx_k\geq 0$.  To see
this, note that the eigenvalues of $\tf 12(H_k+H_k^\top)$ are the same
as the real parts of the eigenvalues of $H_k$. The eigenvectors and
eigenvalues of $H_k$ can be computed as follows.  Let $\omega$ be the
cube root of 1 that is $\omega = -\tf 12+i\tf 12 \sqrt{3}$. Then
\begin{align*}
V_k&=V_1\otimes V_{k-1}\,,\quad
\text{where\ }
V_1=\begin{pmatrix}1&1&1\\1&\omega&\omega^2\\1&\omega^2&\omega\end{pmatrix}\,.
\end{align*}
We write $V_k=U_k+iW_k$, and shall make use of the facts that
$U_k=U_1\otimes U_{k-1}-W_1\otimes W_{k-1}$ and $W_k=U_1\otimes
W_{k-1}+W_1\otimes U_{k-1}$.  It is easily checked that the
eigenvectors of $P_i$ are the columns (and rows) of the symmetric
matrix $V_k$ and that the first row of $V_k$ is
$(1,1,\ldots,1)$. The eigenvalues are also supplied in $V_k$, because
if $V_k(j)$ denotes the $j$th column of $V_k$ (an eigenvector), we
have $P_iV_k(j)=V_k(i,j)V_k(j)$. Thus the corresponding eigenvalue is
$V_k(i,j)$.  Since $H_k$ is a sum of the $P_i$, we also have
$H_kV_k(j)=\sum_ix_iV_k(i,j)V_k(j)$, so the eigenvalue is
$\sum_ix_iV_k(i,j)$, or $\sum_i V_k(j,i)x_i$ since $V_k$ is symmetric.
Thus the real parts of the eigenvalues of $H_k$ are the elements of
the vector $U_kx_k$. This is nonnegative if and only if the symmetric
matrix $\tf 12(H_k+H_k^\top)$ is positive semidefinite.
\smallskip

Recall that $1_k$ denotes the length $3^k$ vector of 1s. 
We will show that we may take $H_k=\sum_i x_k(i)P_i$, where
\begin{align}
x_1 & = (2,2,1)^\top\nonumber\\
x_2 & = (3, 3, 2, 3, 3, 2, 1, 1, 0)^\top\nonumber\\
\intertext{and that we may choose $a_k$ so that for $k\geq 3$,}
x_{k}&=1_{k}+(1,0,0)^\top\otimes x_{k-1}
+(0,1,0)^\top\otimes(a_k,a_k,2,2,a_k,2,1,1,1)^\top\otimes 1_{k-3}\,.\label{xcalc}
\end{align}
In this construction of $x_k$ the parameter $a_{k}$ is chosen
maximally such that $U_kx_k\geq 0$ and $m_k\geq x_k$.\footnote{There
are many choices of $x_2$ that work. We can also take $x_2=(3,3,2,
2,3,2,1,1,1)$ or $x_2=(3,3,2,3,2,2,1,1,1)$.}  The sum of the components of
$x_{k}$ is
\[
1_{k}^\top x_{k} = 3^{k}+1_{k-1}^\top x_{k-1}+3^{k-2}(3+a_k)\,.
\]
To prove the theorem we want $1_k^\top x_k/3^k=w_k$, where these are
the values specified in \eqref{wkcalc}. This requires the values of
the $a_k$ to be:
\begin{equation}
a_k = \left\{\begin{array}{ll}
3-\dfrac{1}{3^{(k-3)/2}}\,,          & \text{when $k$ is odd,}\\[6pt]
3-\dfrac{2}{3^{(k-2)/2}}\,,          & \text{when $k$ is even.}
	     \end{array}\right.\label{akcalc}
\end{equation}\smallskip
So
\[
\{a_3,a_4,\ldots,a_{11},\ldots\}=\{2,\tf{7}{3},\tf{8}{3},
  \tf{25}{9},\tf{26}{9},
  \tf{79}{27},\tf{80}{27},
  \tf{241}{81},\tf{242}{81},\ldots\}\,.
\]
Alternatively, the values of $3-a_k$ are $1,\tf 23,\tf 13,\tf 29,\tf
19,\tf 2{27},\ldots\ $.  For example, with $a_3=2$ we have
\begin{align*}
m_3&=(4,4,3,4,4,3,2,2,2,4,4,3,4,4,3,2,2,2,1,1,1,1,1,1,1,1,1)\,,\\
x_3&=(4,4,3,4,4,3,2,2,1,3,3,3,3,3,3,2,2,2,1,1,1,1,1,1,1,1,1)\,.
\end{align*}
Note that $a_k$ increases monotonically in $k$, from 2 towards 3.  As
$k\rightarrow\infty$ we find $a_k\rightarrow 3$ and $1_{k}^\top
x_{k}/3^k\rightarrow \tf 52$. It remains to prove that with these $a_k$
we have $m_k\geq x_k$ and $U_kx_k\geq 0$.\smallskip

\subsection*{\boldmath $m_k\geq x_k$} 
To prove $m_k\geq x_k$ is easy; we use induction. The base of the
induction is\\ $m_2=(3, 3, 2, 3, 3, 2, 1, 1, 1)\geq
x_2=(3,3,2,3,3,2,1,1,0)$. Assuming $m_{k-1}\geq x_{k-1}$, we then have
\begin{align*}
m_{k}&=1_{k}+(1,1,0)^\top \otimes m_{k-1}\\
&\geq 1_{k}+(1,0,0)^\top\otimes
m_{k-1}+(0,1,0)^\top\otimes\bigl(1_{k-1}\\
&\quad\quad\quad+(1,1,0)^\top\otimes 1_{k-2}
+(1,1,0,1,1,0,0,0,0)^\top\otimes 1_{k-3}\bigr)\\
&=1_{k}+(1,0,0)^\top\otimes m_{k-1}+(0,1,0)^\top\otimes
(3,3,2,3,3,2,1,1,1)^\top\otimes 1_{k-3}\\
&\geq 1_{k}+(1,0,0)^\top\otimes x_{k-1}\
+(0,1,0)^\top\otimes(a_k,a_k,2,2,a_k,2,1,1,1)^\top\otimes 1_{k-3}\\
&= x_{k}\,.
\end{align*}

\subsection*{\boldmath $U_kx_k\geq 0$ } 
To prove $U_kx_k\geq 0$ is much harder. Indeed, $U_kx_k$ is barely
nonnegative, in the sense that as $k\rightarrow\infty$, $\tf 5 9$ of
its components are $0$, and $\tf 29$ of them are equal to $\tf
32$. Thus most of the eigenvalues of $\tf 12(H_k+H_k^\top)$ are $0$. We
do not need this fact, but it is interesting that $2\,U_kx_k$ is a
vector only of integers.

Let $f_k$ be a vector of length $3^k$ in which the first component is
1 and all other components are 0.  Using the facts that
$U_{k}=U_1\otimes U_{k-1}- W_1\otimes W_{k-1}=U_3\otimes U_{k-3}-
W_3\otimes W_{k-3}$ and $W_k1_k=0$ and $U_k1_k=3^kf_k$, we have
\begin{align}
U_2x_2 &= (18,\tf 32,\tf 32,3,0,0,3,0,0)^\top\,,\nonumber\\
\intertext{and for $k\geq 3$,}
U_{k}x_{k} &=3^{k}f_{k}+(1,1,1)^\top\otimes U_{k-1}x_{k-1}\nonumber\\
&+\Bigl(U_3\,\bigl(
(0,1,0)^\top\otimes(a_{k},a_{k},2,2,a_{k},2,1,1,1)^\top\bigr)
\Bigr)\otimes
U_{k-3}1_{k-3}\nonumber\\[3pt]
&=3^{k}f_{k}+(1,1,1)^\top\otimes U_{k-1}x_{k-1}
 +3^{k-3}r_k\otimes f_{k-3}\,,\label{ucalc}
\end{align}
where $r_k$ is
\begin{align}
r_k&=U_3\,\bigl(
(0,1,0)\otimes(a_{k},a_{k},2,2,a_{k},2,1,1,1)\bigr)^\top\nonumber\\[3pt]
&=
\tf 32 \bigl(6+2a_k,0,0,a_k-1,0,a_k-2,a_k-1,a_k-2,0,\label{r1}\\[3pt]
&\quad\quad\quad\quad\quad\quad\quad  -3-a_k,2-a_k,a_k-2,-a_k,0,0,1,2-a_k,0,\label{r2}\\
&\quad\quad\quad\quad\quad\quad\quad\quad\quad\quad\quad\quad\quad\quad -3-a_k,a_k-2,2-a_k,1,0,2-a_k,-a_k,0,0\bigr)^\top\,.\label{rcalc}
\end{align}
Note that we make a small departure from our subscripting convention,
since $r_k$ is not of length $3^k$, but of length $27$. We use the
subscript $k$ to denote that $r_k$ is a function of $a_k$.\smallskip

Using \eqref{ucalc}--\eqref{rcalc} it easy to compute the values
$U_kx_k$, for $k=2,3,\ldots\ $. Notice that there is no need to
calculate the $3^k\times 3^k$ matrix $U_k$. Computing $U_kx_k$ as far
as $k=15$, we find that for the values of $a_k$ conjectured in
\eqref{akcalc} we do indeed always have $U_{k}x_{k}\geq 0$. This gives a
lower bound on the rendezvous value of $w\geq w_{15}=16400/6561\approx
2.49962$. It would not be hard to continue to even larger $k$
(although $U_{15}x_{15}$ is already a vector of length
$3^{15}=14,348,907$). Clearly the method is working. It now remains to
prove that $U_{k}x_{k}\geq 0$ for all $k$.\smallskip

Consider the first third of $U_kx_k$. This is found from 
\eqref{xcalc} and \eqref{r1} to be
\[
3^k f_{k-1}+U_{k-1}x_{k-1}+3^{k-3}\tf 32
\bigl(6+2a_k,0,0,a_k-1,0,a_k-2,a_k-1,a_k-2,0\bigr)\otimes f_{k-3}\,.
\]
Assuming $U_{k-1}x_{k-1}\geq 0$ as an inductive hypothesis, and using
the fact that $a_k\geq 2$, this vector is nonnegative. So this part of 
$U_kx_k$ is nonnegative.

As for the rest of $U_kx_k$ (the part that can be found from \eqref{xcalc} 
and \eqref{r2}--\eqref{rcalc}), notice that $r_k$ is symmetric, in the sense
that $S_3r_k=r_k$, where
\[
S_1=\begin{pmatrix}1&0&0\\0&0&1\\0&1&0\end{pmatrix} 
\]
and $S_3=S_1\otimes S_1\otimes S_1$.  The matrix $S_k$ transposes 1s
and 2s. Indeed $S_kP_i=P_i^\top$. Thus the proof is complete if we can
show that just the middle third of $U_kx_k$ is nonnegative.  Assuming
that $U_{k-1}x_{k-1}\geq 0$ and $a_k\geq 2$, there are just 4
components of this middle third that depend on $a_k$ and which might
be negative.  Let $I_{k}$ denote a $3^k\times 3^k$ identity matrix.
This middle third is found from \eqref{xcalc} and \eqref{r2} and is as
follows, where we indicate in bold face terms that might be negative,
\begin{align*}
\bigl((0,1,0)\otimes &I_{k-1}\bigr)U_kx_k\\
&=U_{k-1}x_{k-1}+\tf 32 3^{k-3}\bigl( \bm{-3-a_k},\bm{2-a_k},a_k-2,\bm{-a_k},0,0,1,\bm{2-a_k},0
\bigr)\otimes f_{k-3}\,.
\end{align*}
The four possibly negative components of the middle third are shown
above in bold and are
\begin{align}
t_{k1}&=(0,1,0)\otimes(1,0,0,0,0,0,0,0,0)\otimes f_{k-3}^\top\,U_kx_k\nonumber\\
&=(U_{k-1}x_{k-1})_1+\tf 32 3^{k-3}\left(-3- a_k\right)\label{t1}\\[5pt]
t_{k2}&=(0,1,0)\otimes(0,1,0,0,0,0,0,0,0)\otimes f_{k-3}^\top\,U_kx_k\nonumber\\
&=(U_{k-1}x_{k-1})_{3^{k-3}+1}+\tf 32 3^{k-3}\left(2-a_k\right)\label{t2}\\[5pt]
t_{k3}&=(0,1,0)\otimes(0,0,0,1,0,0,0,0,0)\otimes f_{k-3}^\top\,U_kx_k\nonumber\\
&=(U_{k-1}x_{k-1})_{3\, 3^{k-3}+1}+\tf 32 3^{k-3}(-a_k)\label{t3}\\[5pt]
t_{k4}&=(0,1,0)\otimes(0,0,0,0,0,0,0,1,0)\otimes f_{k-3}^\top\,U_kx_k\nonumber\\
&=(U_{k-1}x_{k-1})_{7\, 3^{k-3}+1}+\tf 32 3^{k-3}\left(2-a_k\right)\label{t4}
\end{align}
The remainder of the proof is devoted to proving that all these are nonnegative.
Consider $t_{k1}$.  It is easy to work out a formula for $t_{k1}$,
since
\begin{align*}
(U_kx_k)_1 &=f_k^\top U_kx_k\\
&=3^k+f_{k-1}^\top U_{k-1}x_{k-1}+3^{k-3}\tf 32(6+2a_k)\\[4pt]
&= (U_{k-1}x_{k-1})_1 +4\, 3^{k-1}+3^{k-2}a_k
\end{align*}
Thus
\begin{equation}
(U_kx_k)_1 = 2\, 3^k+\sum_{i=3}^k 3^{i-2}a_i\,,\label{fcalc}
\end{equation}
and
\begin{equation}
t_{k1}=\tf 12\, 3^{k}+\sum_{i=3}^{k-1} 3^{i-2}a_i-\tf 12 3^{k-2}a_k
\end{equation}
This is  nonnegative since $a_k\leq 3$.

Amongst the remaining terms, we observe empirically that $t_{k2}\geq
t_{k4}\geq t_{k3}$.  It is $t_{k3}$ that is the least of the four
terms, and which constrains the size of $a_k$. Let us begin therefore by
finding a formula for $t_{k3}$. We have
\begin{align*}
t_{k3}&=(0,1,0)\otimes(0,0,0,1,0,0,0,0,0)\otimes f_{k-3}^\top\,U_kx_k\nonumber\\
&=(0,0,0,1,0,0,0,0,0)\otimes f_{k-3}^\top U_{k-1}x_{k-1}-3^{k-2}\tf 12 a_k\\
&=(0,1,0)\otimes f_1^\top\otimes f_{k-3}^\top\bigl(
3^{k-1}f_{k-1}+(1,1,1)^\top\otimes U_{k-2}x_{k-2}
 +3^{k-4}r_{k-1}\otimes f_{k-4}
\bigr)
-3^{k-2} \tf 12 a_k\\
&=(1,0,0,0,0,0,0,0,0)\otimes f_{k-4}^\top  U_{k-2}x_{k-2}+3^{k-4}(0,1,0)\otimes f_2)r_{k-1}
-3^{k-2}\tf 12 a_k\\
&=(U_{k-2}x_{k-2})_1 - 3^{k-4}\tf 32(3+ a_{k-1})-3^{k-2}\tf 12 a_k\\
&=(U_{k-2}x_{k-2})_1 - 3^{k-3}\tf 12(3+ a_{k-1})-3^{k-2}\tf 12 a_k
\end{align*}
This means that $t_{k3}$ can be computed from the first component of
$U_{k-2}x_{k-2}$, which we have already found in \eqref{fcalc}. So
\begin{align}
t_{k3}&=2\, 3^{k-2}+\sum_{i=3}^{k-2} 3^{i-2}a_i- 3^{k-3}\tf 12(3+
a_{k-1})-3^{k-2}\tf 12 a_k\nonumber\\
&=\tfrac{1}{2}3^{k-1}+\sum_{i=3}^{k-2} 3^{i-2}a_i
-\tf 12 3^{k-3}a_{k-1}-\tf 12 3^{k-2}a_{k}\,.\label{tk3calc}
\end{align}
We now put the $a_k$ to the values specified in \eqref{akcalc}.
It is easy to check with \eqref{akcalc} and \eqref{tk3calc} that
$t_{k3}=0$ for all $k$.
\smallskip

It remains only to check that also $t_{k2}\geq 0$ and $t_{k4}\geq
0$. We have
\begin{align*}
t_{k2}&=(0,1,0)\otimes(0,1,0,0,0,0,0,0,0)\otimes f_{k-3}^\top\,U_kx_k\nonumber\\
&=(0,1,0,0,0,0,0,0,0)\otimes f_{k-3}^\top U_{k-1}x_{k-1}+3^{k-2}(1-\tf 12 a_k)\\
&=(1,0,0)\otimes(0,1,0)\otimes f_{k-3}^\top\bigl(
3^{k-1}f_{k-1}+(1,1,1)^\top\otimes U_{k-2}x_{k-2} 
+3^{k-4}r_{k-1}\otimes f_{k-4}
\bigr)\\
&\quad +3^{k-2}(1-\tf 12 a_k)\\
&=(0,1,0)\otimes f_{k-3}^\top U_{k-2}x_{k-2}-3^{k-4}\tf 32(1-a_{k-1})+3^{k-2}(1-\tf 12 a_k)\,.
\end{align*}
We recognize $(0,1,0)\otimes f_{k-3}^\top U_{k-2}x_{k-2}$ to be the
first component of the middle third of $U_{k-2}x_{k-2}$. The
recurrence relation for this is
\begin{align*}
(0,1,0)\otimes f_{k-1}^\top U_{k}x_{k}
&=(0,1,0)\otimes f_{k-1}^\top\bigl(3^{k}f_{k}+(1,1,1)^\top\otimes
  U_{k-1}x_{k-1}
+3^{k-3}r_k\otimes f_{k-3}\Bigr)\\
&=f_{k-1}^\top U_{k-1}x_{k-1}-3^{k-2}\tf 12(3+a_{k})\,.
\end{align*}
The right hand side can be computed from \eqref{fcalc}.
So we now have,
\begin{align}
t_{k2} &= 2\, 3^{k-3}+\sum_{i=3}^{k-3}3^{i-2} a_i
-3^{k - 4}\tf 12(3 + a_{k-2}) - 3^{k - 3}\tf 12(1 - a_{k - 1}) + 
  3^{k - 2}(1 - \tf 12a_k)\nonumber\\
&=4\, 3^{k-3} +\sum_{i=3}^{k-3}3^{i-2} a_i
-\tf 12 3^{k-4}a_{k-2}
+\tf 12 3^{k-3}a_{k-1}
-\tf 12 3^{k-2}a_{k}\,.\label{tk2calc}\\
\intertext{Finally, we establish a formula for $t_{k4}$.}
t_{k4}&=(0,1,0)\otimes(0,0,0,0,0,0,0,1,0)\otimes f_{k-3}^\top\,U_kx_k\nonumber\\
&=(0,0,0,0,0,0,0,1,0)\otimes f_{k-3}^\top U_{k-1}x_{k-1}+3^{k-2}(1-\tf 12 a_k)\nonumber\\
&=(0,0,1)\otimes(0,1,0)\otimes f_{k-3}^\top\bigl(
3^{k-1}f_{k-1}+(1,1,1)^\top\otimes U_{k-2}x_{k-2}
 +3^{k-4}r_{k-1}\otimes f_{k-4}
\bigr)\\
&\quad +3^{k-2}(1-\tf 12 a_k)\nonumber\\
&=(0,1,0)\otimes f_{k-3}^\top U_{k-2}x_{k-2}+3^{k-4}\tf
32+3^{k-2}(1-\tf 12 a_k)\nonumber\\
&=5\, 3^{k-3} +\sum_{i=3}^{k-3}3^{i-2} a_i
-\tf 12 3^{k-4}a_{k-2}
-\tf 12 3^{k-2}a_k\,.\label{t4calc}
\end{align}
Thus we can check the fact that we observed empirically, that
$t_{k2}\geq t_{k4}\geq t_{k3}$. We find
\begin{align*}
t_{k2}-t_{k4} & =\tf 12 3^{k-3}(a_{k-1}-2)\,,\\[3pt]
t_{k4}-t_{k3} & =\tf 12 3^{k-3}(1-a_{k-2}+a_{k-1})\,.
\end{align*}
Since $a_k$ is at least 2 and $a_k$ is increasing in $k$, both of the
above are nonnegative. So $t_{k2}$ and $t_{k4}$ are both at least as
great as $t_{k3}$, which we have already shown to be $0$.  This
establishes $U_kx_k\geq 0$ and so the proof is now complete.\hfill
\rule{9pt}{9pt}

\section{On discovery of the proof}

The proof begs the question: how did we guess the recursion for $x_k$?
Let us restate it here for convenience. With $a_k$ given by
\eqref{akcalc}, the recursion is
\begin{equation}
x_k = 1_k + (1,0,0)^\top\otimes
x_{k-1}+(0,1,0)\otimes(a_k,a_k,2,2,a_k,2,1,1,1)\otimes 1_{k-3}\,.\label{xcalc2}
\end{equation}
Let us briefly describe the
steps and ideas in research that led to \eqref{xcalc2}. 

We began by computing lower bounds on $w_k$ by solving the
semidefinite programming problem
\begin{equation}
\maximize\ \tr(J_kH_k) \,:\, H_k\leq M_k\,,\ H_k\succeq 0\,.\label{sdp}
\end{equation}
A similar line of approach has been followed concurrently by Han, Du,
Vera and Zuluaga (2006). The lower bounds that are obtained by solving
\eqref{sdp} turn out to be achieved by the \AW strategy and so are
useful in proving the Fan--Weber conjecture (that \AW minimizes
$E[\min\{T,k+1\}]$) up to $k=5$. However, they only produce numerical
answers, with little guide as to a general form of solution. In fact,
since one can only solve the SDPs up to the numerical accuracy of a
SDP solver (which, like \emph{sedumi}, uses interior point methods),
such proofs are only approximate. For example, by this method one can
only prove that $w_5\geq 2.40740740$, but not $w_5=\tf {65}
{27}=2.4\dot{\overline{074}}$.

We tried to find rational solutions so the proofs could be exact.  A
major breakthrough was to realise that we could compute a common
eigenvector set for $P_1,\ldots,P_{3^k-1}$ and write $M_k=\sum_i
m_k(i)P_i$. We discovered this as we noticed and tried to explain the
fact that the real parts of all the eigenvalues of $2\, M_k$ are
integers. (In fact, there is a little-known theorem which says that if
a real symmetric matrix has only integer entries then all its rational
eigenvalues must be integers.)  This allowed us to recast \eqref{sdp}
as the linear program
\begin{equation}
\maximize\ \sum_{i=0}^{3^k-1} x(i) \,:\, x\leq m_k\,,\ U_kx\geq 0\,.\label{lp}
\end{equation}
Now we could find exact proofs of the Fan--Weber conjecture as far as
$k=8$ , where $U_8$ is $6561\times 6561$. These solutions were found
using \emph{Mathematica} and were in rational numbers, thus 
providing us with tight proofs of the optimality of \AW up to $k=8$.
They also allowed us to prove the Fan-Weber conjecture for greater
values of $k$  since the number of decision variables in the LP
grows as $3^{k}$, whereas in the SDP it grows as $3^{2k}$.


It seems very difficult to find a general solution to \eqref{lp} that
will hold for all $k$. The LP is highly degenerate with many optimal
solutions. There are indeed $12$ different extreme point solutions to
the LP at just $k=2$. No general pattern to the solution emerges as it
is solved for progressively larger $k$.  For, say $k=4$, there are
many $H_4$ that can be used to prove $w_k=\tf 73$. We searched amongst
the many solutions for ones with some pattern that might be
generalized.  This proved very difficult. We tried forcing lots of
components of the solution vector to be integers, or identical, and
looked for solutions in which the solution vector for $k-1$ was
embedded within the solution vector for $k$. We looked at adding other
constraints, and constructed some solutions by augmenting the
objective function and choosing amongst possible solution by a
minimizing a sum of squares penalty.

Another approach to the problem of minimizing $p^\top M_kp$ over
$p\in\Delta_k$ is to make the identification $Y=pp^\top$.  With this
identification, $Y$ is positive semidefinite, $\tr(J_kY)=1$, and
$\tr(M_kY)= \tr(M_kpp^\top)=p^\top M_k p$. This motivates a
semidefinite programming relaxation of our problem: minimize
$\tr(M_kY)$, subject to $\tr(J_kY)=1$ and $Y\succeq 0$. This can be
recast as the LP
\begin{equation}
\minimize\ y^\top m_k \,:\, y^\top U_k\geq 0\,,\ 1_k^\top y=1\,,\ y\geq
0\,. \label{duallp}
\end{equation}
This is nearly the dual of \eqref{lp} (which is the same, but has an
additional constraint of $y^\top U_kS_k=y^\top U_k$).

With \eqref{duallp} in mind, we imagined taking $y$ as \AW and worked
at trying to guess a full basis in the columns of $U_k$ that is
complementary slack to $y$ and from which one can then compute a
solution to \eqref{lp}.  We also explored a number of different LP
formulations.  All of this was helpful in building up intuition as to
how a general solution might possibly be constructed.

Another major breakthrough was to choose to work with the constraint
$x\leq m_k$ in which $m_k$ is the first row of the nonsymmetric matrix
$M_k$, rather than to use the first row of the symmetric matrix $\tf
12 (M_k+M_k^\top)$. By not `symmetrizing' $M_k$ we were able to find
solutions with a simpler form, and felt that there was more hope in
being able to write the solution vector $x_k$ in a Kronecker product
calculation with the solution vector $x_{k-1}$. Noticing that all the
entries in $M_k$ are integers, we found that it was possible to find a
solution for $H_k$ in which all the entries in $H_k$ are integers, as
far as $k=5$. It is not known whether this might be possible for even
greater $k$. The $H_k$ constructed in the proof above have entries
that are not integers, although they are of course rational.

Since $M_k$ is computed by Kronecker products it is natural to look
for a solution vector of a form in which $x_k$ is expressed in terms
of $x_{k-1}$ in some sort formula using Kronecker products. The final
breakthrough came in discovering the length 27 vector
$(0,1,0)\otimes(a_k,a_k,2,2,a_k,2,1,1,1)$. This was found only after
despairing of something simpler. We expected that if it were possible
to find a Kronecker product form solution similar to \eqref{xcalc2},
then this would use a vector like the above, but of length only 3 or
9. It was only when we hazarded to try something of length 27 that the
final pieces fell in place. The final trick was to make the formula
for obtaining $x_k$ from $x_{k-1}$ not be constant, but depending on
$k$, as we have done with our $a_k$. We were lucky at the end that we
could solve the recurrence relations for $t_{k1},t_{k2},t_{k3},t_{k4}$
and prove $U_kx_k\geq 0$. It all looks so easy with hindsight!

\section{Ongoing research}

\begin{enumerate}
\item One would like to have a direct proof that $w=\tf 52$, without
needing to also find the $w_k$. Perhaps an idea for such a proof is
pregnant within the proof above.\medskip

\item While for many graphs it is possible to use the solution of
a semidefinite programming problem to obtain a lower bound on the
rendezvous value, it is not usually possible to recast the SDP as a
linear program. A very important feature of the $K_3$ problem is that
it is so strongly captured within the algebra of a 
group of rotational symmetry, whose permutation matrices are the
$P_i$. This continues to be true for rendezvous search on $C_n$, in
which $n$ locations are arranged around a circle and players have a
common notion of clockwise. We are presently looking for results in
that direction.


\item It is as easy consequence of Theorem 1 that \AW maximizes $E[\beta^T]$
for all $\beta\in(0,1)$. This follows from the fact that \AW minimizes
$\sum_{i=0}^k P(T>i)$ for all $k$.

\item It will be interesting to explore whether our methods are
  helpful for rendezvous problems on $K_n$ ($n\geq 4$), or on other
  graphs. It is not hard to compute the optimal \AW strategy for
  $K_n$. See Anderson and Weber (1990). For example, for $n=4$, an \AW
  strategy achieves $ET\approx 3.4247$, using probabilities of
  staying and touring of $p\approx 0.3320$ and $1-p\approx 0.6680$,
  respectively. As $n\rightarrow\infty$, an \AW strategy achieves
  $ET \approx .8289\,n$ with $p\approx 0.2475$. Interestingly, Fan
  (2009) has shown that if the rendezvous game is played on $K_4$,
  locations are imagined to be placed around a circle, and players are
  provided with a common notion of clockwise, then there exists a
  $3$-Markovian strategy that is better than \AW. However, it is
  open as to whether \AW is optimal when players are
  not provided with such initial information.

\item We conjecture that \AW is optimal in a rendezvous game played on
  $K_3$ in which players may overlook one another with probability
  $\epsilon$, (that is, they can fail to meet even when they are in
  the same location). This is easily shown to be true for the game on
  $K_2$. To analyse this game on $K_3$ we simply redefine
\[
B_1=
\begin{pmatrix}
1 & 1 & \epsilon\\
\epsilon & 1 & 1\\
1 & \epsilon & 1
\end{pmatrix},
\]
where $0<\epsilon<1$. We can generalize all the ideas in this paper,
except that we have not been able to guess a construction for the
matrix $H_k$. Fan (2009) has observed that not only does \AW appear to
be optimal, but also that the optimal probability of `staying' is the
same for all $\epsilon$, i.e., $p=1/3$. However, for games on $K_4,
K_5,\ldots,\,$ the optimal value of $p$ is
decreasing in $\epsilon$.
\medskip

\item In the asymmetric version of the rendezvous search game (in
 which players I and II can adopt different strategies) the rendezvous
 values for the games on $K_2$ and $K_3$ are 1 and 1.5 respectively
 (and are achieved by the `wait-for-mommy' strategy). These are
 exactly 1 less than the rendezvous values of 2 and 2.5 that pertain
 in the symmetric games (and are achieved by the \AW strategy).

 A rendezvous search game can also be played on a line. The players
 start 2 units apart and can move 1 unit left or right at each
 step. The asymmetric rendezvous value is known to be $3.25$ (Alpern
 and Gal, 1995). In the symmetric game it is known that $4.1820\leq
 w\leq 4.2574$.  Han, et al. (2006) have conjectured $w=4.25$. If that
 is correct then the difference in rendezvous values for asymmetric and
 symmetric games is again exactly 1.

\item In the symmetric rendezvous search game played on $K_3$ it is of
no help to the players to be provided with a common notion of
clockwise. Similarly, our experience in studying the symmetric
rendezvous search game on the line suggests that it is no help to the
players to be provided with a common notion of left and right.

\item No one has yet found a way to prove that the rendezvous value for
the symmetric rendezvous search game on $K_n$ is an increasing
function of $n$.
\end{enumerate}

\subsection*{Thanks} 

I warmly thank my Ph.D.\ student Junjie Fan (Jimmy) for his
enthusiasm, many helpful discussions and proof-reading of this paper.
By pioneering the use of semidefinite programming as a method of
addressing rendezvous search problems, he has been the first in many
years to obtain significantly improved lower bounds on $w$.  



\end{document}